\documentclass[11pt]{article}
 \input xy
 \xyoption{all}
 \newcommand{\ver}[1]{*+<2mm>[F-:<3pt>]{#1}}

 \usepackage{esint}

 \usepackage{fullpage}

\newcommand{\rvsR}{rank-varying sub-Riemannian structure}
\newcommand{\ov}{\overline}

\newcommand{\ffoot}[1]{}

\newcommand{\rotations}{revolutions}

\newcommand{\HH}{{\bf (H0)}}
\newcommand{\same}{{equivalent}}
\newcommand{\PH}{\varphi}
\newcommand{\bD}{\Delta}
\newcommand{\pphi}{\xi}

\renewcommand{\th}{\theta}

\newcommand{\cil}{{\bf W}}

\newcommand{\f}{f}

\usepackage{mathrsfs}
\usepackage{graphicx}
\usepackage{amssymb}
\usepackage{color}
\newtheorem{theorem}{Theorem}
\newtheorem{corollary}[theorem]{Corollary}
\newtheorem{lemma}[theorem]{Lemma}
\newtheorem{proposition}[theorem]{Proposition}
\newtheorem{definition}[theorem]{Definition}
\newtheorem{remark}[theorem]{Remark}

\newcommand{\bt}{\begin{theorem}}
\newcommand{\et}{\end{theorem}}
\newcommand{\bl}{\begin{lemma}}
\newcommand{\el}{\end{lemma}}
\newcommand{\bp}{\begin{proposition}}
\newcommand{\ep}{\end{proposition}}
\newcommand{\bc}{\begin{corollary}}
\newcommand{\ec}{\end{corollary}}
\newcommand{\bdeff}{\begin{definition}}
\newcommand{\edeff}{\end{definition}}
\newcommand{\brem}{\begin{remark}}
\newcommand{\erem}{\end{remark}}
\renewcommand{\r}[1]{(\ref{#1})}
\newcommand{\con}{{\mathcal C}}

\newcommand{\bi}{\begin{itemize}}
\newcommand{\iii}{\item}
\newcommand{\ei}{\end{itemize}}
\newcommand{\bd}{\begin{description}}
\newcommand{\ed}{\end{description}}

\newcommand{\bqn}{\begin{eqnarray}}
\newcommand{\eqn}{\end{eqnarray}}
\newcommand{\eqnn}{\nonumber\end{eqnarray}}

\newcommand{\nn}{\nonumber}
\newcommand{\ba}[1]{\begin{array}{#1}}
\newcommand{\ea}{\end{array}}

\newcommand{\R}{\mathbb{R}}

\newcommand{\g}{\gamma}

\newcommand{\rvd}{rank-varying distribution}

\newcommand{\VecM}{\mathrm{Vec}(M)}

\newcommand{\Gq}{{\gg}_q}% cerca le altre

\newcommand{\CC}{{\mathfrak C}}
\newcommand{\Gg}{{\cal G}}
\newcommand{\V}{{\cal V}}
\newcommand{\Ee}{{\cal E}}

\newcommand{\Id}{{\mbox Id}}

\newcommand{\Zz}{\mathcal{Z}}
\renewcommand{\gg}{{\bf G}}
\newcommand{\sign}{\mathrm{sign}}

\title{\LARGE \bf
Lipschitz Classification of  Almost-Riemannian Distances on Compact Oriented Surfaces
}

\author{ 
 U.~Boscain\thanks{
Centre de Math\'ematiques Appliqu\'ees, \'Ecole Polytechnique
Route de Saclay, 91128 Palaiseau Cedex, France,
        {\tt ugo.boscain@cmap.polytechnique.fr} The author has been supported by the ERC Starting Grants GeCoMethods, contract 239748}
, G.~Charlot\thanks{
Institut Fourier, UMR 5582, CNRS/Universit\'e Grenoble 1, 100 rue des Maths, BP 74, 38402 St Martin d'H\`eres, France,
{\tt Gregoire.Charlot@ujf-grenoble.fr}} , 
R.~Ghezzi\thanks{
 SISSA, via Beirut 2-4, 34014 Trieste, Italy, {\tt ghezzi@sissa.it}}, M.~Sigalotti\thanks{INRIA Nancy -- Grand Est, \'Equipe-projet CORIDA, and
Institut \'Elie Cartan, UMR CNRS/INRIA/Nancy Universit\'e, BP 239, 54506 Vand\oe uvre-l\`es-Nancy, France,
        {\tt Mario.Sigalotti@inria.fr}} 
}

\begin{document}
\maketitle
\begin{abstract}
Two-dimensional almost-Riemannian structures are generalized Riemannian structures on surfaces for which a local orthonormal frame is given by a  Lie bracket generating pair of vector fields that can become collinear.
We consider the  Carnot--Caratheodory distance canonically associated with an almost-Riemannian
 structure and study the problem of Lipschitz equivalence between two such distances on the same compact oriented surface.
 We analyse the generic case, allowing in particular for the  presence of tangency points, i.e., points where two generators of the distribution and their Lie bracket are linearly dependent.  The main result of the paper provides a characterization of the Lipschitz equivalence class  
of an almost-Riemannian distance in terms of a labelled graph associated with it.
\end{abstract}

\section{Introduction}

Consider a pair of smooth vector fields $X$ and $Y$ on a two-dimensional smooth manifold $M$.
If the pair $(X,Y)$ is Lie bracket generating, i.e., if
$\mathrm{span}\{X(q),Y(q),[X,Y](q),[X,[X,Y]](q),\ldots\}$ is
full-dimensional at every $q\in M$, then the
control system
\bqn
\label{ff}
\dot q=u X(q)+v Y(q)\,,~~~u^2+v^2\leq 1\,,~~~q\in M\,,
\eqn
is completely controllable and  the
minimum-time function defines a
continuous distance $d$ on $M$.
When $X$ and $Y$ are everywhere linear independent
(the only possibility for this to happen is that $M$ is parallelizable),
such distance is Riemannian and it
corresponds
to the metric for which $(X,Y)$ is an orthonormal  frame.
Our aim is to study the geometry obtained starting from a pair of vector fields which
may  become collinear.
Under generic hypotheses, the set $\Zz$ (called {\it singular locus}) of
points of $M$ at which  $X$ and $Y$ are parallel
is  a one-dimensional embedded submanifold of $M$ (possibly disconnected).

Metric structures
that can be defined {\it locally} by a pair of vector fields
$(X,Y)$ through \r{ff}
are called almost-Riemannian structures. 

Equivalently, an almost-Riemannian structure ${\cal S}$  can be defined  as an Euclidean  bundle $E$ of rank two over $M$ (i.e. a vector bundle whose fibre is equipped with a smoothly-varying scalar product $\langle\cdot,\cdot\rangle_q$) and a morphism of vector bundles $\f:E\rightarrow TM$ such that the evaluation at $q$ of the Lie algebra generated by the submodule 
\bqn
\Delta:=\{\f\circ\sigma\mid\sigma\mathrm{ \, section\,  of\, } E\}\label{delta}
\eqn
of the algebra of vector fields on $M$
 is equal to $T_{q}M$ for every $q\in M$.

 If $E$ is orientable, we say that ${\mathcal S}$ is {\it orientable}. 
The singular locus $\Zz$  is the set of points $q$  of $M$ at which $\f(E_{q})$ is one-dimensional.
An almost-Riemannian structure is Riemannian if and only if $\Zz=\emptyset$, 
i.e. $\f$ is an isomorphism of vector bundles.

 The first example of
genuinely almost-Riemannian structure
is provided by
the Grushin plane,
which is the almost-Riemannian structure on $M=\R^2$
with $E=\R^2\times\R^2$, $\f((x,y),(a,b))=((x,y),(a,b x))$ and $\langle\cdot,\cdot\rangle$ the canonical Euclidean structure on $\R^2$.
The model was originally introduced in the context of
hypoelliptic operator theory
\cite{FL1,grusin1} (see also \cite{bellaiche,algeria}).
Notice that the
singular locus
is indeed nonempty, being equal to
the $y$-axis.
Another example of almost-Riemannian structure
 appeared in problems of control of quantum mechanical systems
(see \cite{q4,q1}). 

Almost-Riemannian structures present very interesting phenomena. For instance,  even in the case where
the Gaussian curvature is everywhere negative (where it is defined, i.e., on $M\setminus\Zz$), geodesics may have conjugate points. 
This happens for instance on the Grushin plane (see \cite{ABS} and also \cite{tannaka,rigge} in the case of surfaces of revolution). The structure of the cut and conjugate loci  is described in \cite{bcgj} under generic assumptions.

In \cite{euler}, we provided an extension of the Gauss--Bonnet theorem to almost-Riemannian structures, linking the Euler number of the vector bundle $E$ to 
a suitable principal part of the integral of the curvature on $M$. For generalizations of the Gauss-Bonnet formula in related context see also \cite{pelletier}.

The results in \cite{euler} have been obtained under a set of  generic hypotheses called \HH. To introduce it, let us define the {\it flag} of the  submodule $\bD$ defined in \r{delta} as the sequence of submodules
$\bD=\bD_1\subset \bD_2\subset\cdots \subset\bD_m \subset \cdots$ 
defined through the recursive formula
$$%\bqn\label{flag}
\bD_{k+1}=\bD_k+[\bD,\bD_k].
$$%\eqn 
Under generic assumptions, the singular locus $\Zz$ has the following properties: {\bf (i)} $\Zz$ is an
embedded one-dimensional 
submanifold of
$M$;
{\bf (ii)} the points $q\in M$ at which $\bD_2(q)$ is
one-dimensional are isolated;
{\bf (iii)}  $\bD_3(q)=T_qM$ for every $q\in M$.
We  say that $\cal S$ satisfies \HH\ if properties {\bf (i)},{\bf (ii)},{\bf (iii)} hold true. If
 this is the case, a point $q$ of $M$ is called {\it ordinary} if $\bD(q)=T_qM$, {\it Grushin point} if $\bD(q)$ is one-dimensional and $\bD_2(q)=T_qM$, i.e. the distribution is transversal to $\Zz$, and {\it tangency point} if $\bD_2(q)$ is one-dimensional, i.e. the distribution is tangent to $\Zz$. Local normal forms around ordinary, Grushin and tangency points have been provided in \cite{ABS}. When an  ARS ${\cal S}=(E,f,\langle\cdot,\cdot\rangle)$ satisfying \HH\ is oriented and the surface itself is oriented, $M$ is split into two open sets $M^+$, $M^-$ such that $\Zz=\partial M^+=\partial M^-$, $f:E|_{M^+}\rightarrow TM^+$ is an orientation preserving isomorphism and $f:E|_{M^-}\rightarrow TM^-$ is an orientation reversing isomorphism. Moreover, in this case it is possible to associate with each tangency point $q$ an integer $\tau_q$ in the following way. Choosing on $\Zz$ the orientation induced by $M^+$,  $\tau_q= 1$ if walking along the oriented curve $\Zz$ in a neighborhood of $q$ the angle between the distribution and the tangent space to $\Zz$ increases, $\tau_q= -1$ if the angle decreases.

\begin{figure}[h!]
\begin{center}
\input{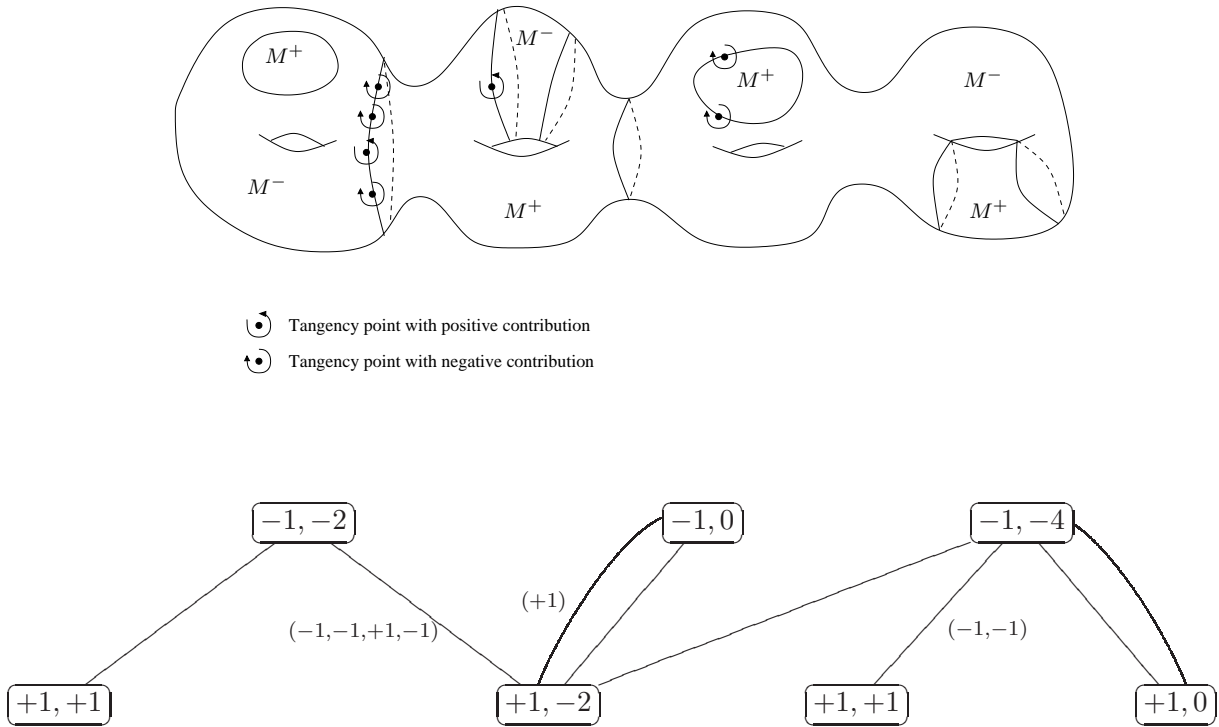}
\vspace{1.5cm}

\xymatrix{
&&\ver{-1,-2}\ar@{-}[ddll]\ar@{-}[ddrr]_{(-1,-1,+1,-1)}&&&\ver{-1, 0}\ar@{-}@(ul,dl)[ddl]_{(+1)}\ar@{-}[ddl]&&\ver{-1,-4}\ar@{-}[ddlll]\ar@{-}[ddl]^{(-1,-1)}\ar@{-}[ddr]\ar@{-}@(ur,dr)[ddr]&\\
&&&&&&\\
\ver{+1,+1}&&&&\ver{+1,-2}&&\ver{+1,+1}&&\ver{+1, 0}
}

\caption{\label{figintro}Example of ARS on a surface of genus $4$ and corresponding labelled graph}
\end{center}
\end{figure}

In this paper we provide a classification of orientable two-dimensional almost-Riemannian structures in terms of graphs. 
With an oriented almost-Riemannian structure, we associate a graph whose vertices correspond to 
connected components of $M\setminus\Zz$ and whose edges correspond to connected components of $\Zz$. 
The edge corresponding to a connected component $W$ of $\Zz$  joins the two vertices corresponding to the connected components of $M\setminus\Zz$ adjacent to $W$. 
Every vertex is labelled with its orientation ($\pm 1$ if it a subset of $M^\pm$) and its Euler characteristic.
Every edge is labelled with the ordered sequence of signs (modulo cyclic permutations) given by the contributions at the tangency points belonging to $W$.  See Figure \ref{figintro} for an example of almost-Riemannian structure and its corresponding graph. %(figure \ref{figintro}(a)). 
We say that two labelled graphs  are equivalent if they are equal or they can be obtained by the same almost-Riemannian structure reversing the orientation of the vector bundle.  

The main result of the paper is the following.
\bt\label{lip-eq}
Two oriented almost-Riemannian structures, defined on compact oriented surfaces and satisfying \HH, are Lipschitz equivalent if and only if they have equivalent graphs.
\et
In the statement above, two  almost-Riemannian structures are said to be Lipschitz equivalent if there exists a diffeomorphism between their base surfaces which is bi-Lipschitz with respect to the two almost-Riemannian distances. 

This theorem shows another interesting difference between  Riemannian manifolds and almost-Riemannian ones: in the Riemannian context, Lipschitz equivalence coincides with the equivalence as differentiable manifolds;  in the almost-Riemannian context, Lipschitz equivalence is a stronger condition. Notice, however, that in general Liptschitz equivalence does not imply isometry. Indeed, the Lipschitz equivalence between two structures does not depend on the metric structure but only on the submodule $\Delta$. This is highlightened by the fact that the graph itself depends only on $\Delta$.

The structure of the paper is the following. In section~\ref{basdef}, we recall some basic
 notion of sub-Riemannian geometry.
Section~\ref{deftau} introduces the definitions of the number of revolution of a one-dimensional distribution along
a closed oriented curve and of the graph associated with an almost-Riemannian structure. In section~\ref{proof} we demonstrate Theorem~\ref{lip-eq}. Section~\ref{necessity} provides the proof of the fact that having equivalent graphs is a necessary condition for Lipschitz equivalent structures. Finally, in section~\ref{sufficiency} we show this condition to be sufficient.

\section{Preliminaries}\label{basdef}

This section is devoted to recall some basic definitions  in the framework of  sub-Riemannian geometry following  \cite{euler,ABS}, see also \cite{bellaiche,montgomery}.

Let $M$ be a $n$-dimensional manifold. Throughout the paper, unless specified,  manifolds are smooth (i.e., $\con^{\infty}$) and without boundary; vector fields  and differential forms are smooth. Given a vector bundle $E$ over $M$, the $\con^\infty(M)$-module of smooth sections of $E$ is denoted by $\Gamma(E)$. For the particular case  $E=TM$, the set of smooth vector fields on $M$ is denoted by $\VecM$. 
\begin{definition}\label{fiberrvd}
A {\it $(n,k)$-\rvd} on a $n$-dimensional manifold $M$ is 
a pair $(E,\f)$ where $E$ is a vector bundle of rank $k$ over $M$ and $\f:E\rightarrow TM$ is a morphism of vector bundles, 
i.e., {\bf (i)}  the diagram 
$$
\xymatrix{
 E  \ar[r]^{\f} \ar[dr]_{\pi_E}   & TM \ar[d]^{\pi}            \\
 & M                          
}    
$$
 commutes, where  $\pi:TM\rightarrow M$ and $\pi_E:E\rightarrow M$ denote the canonical projections and {\bf (ii)} $\f$ is linear on fibers.
Moreover, we require the map $\sigma\mapsto\f\circ\sigma$ from $\Gamma(E)$ to $\VecM$ to be injective.
\end{definition}

Given a $(n,k)$-\rvd, we denote by $f_*:\Gamma(E)\rightarrow \VecM$ the morphism of $\con^\infty(M)$-modules that maps $\sigma\in\Gamma(E)$ to $\f\circ\sigma\in\VecM$. The following proposition shows that all the information about a \rvd\ is carried by the submodule $\f_*(\Gamma(E))$. 
\begin{proposition}\label{stessomodulo} Given two $(n,k)$-\rvd s $(E_i,\f_i), i=1,2$, assume that they define the same submodule of $\VecM$, i.e., $(\f_1)_*(\Gamma(E_1))=(\f_2)_*(\Gamma(E_2))=\bD\subseteq\VecM$. Then, there exists an isomorphism of vector bundles $h:E_1\rightarrow E_2$ such that $\f_2\circ h=\f_1$. 
\end{proposition}
{\bf Proof.} Since $(\f_i)_*:\Gamma(E_i)\rightarrow \bD, i=1,2$, are  isomorphisms of $\con^\infty(M)$-modules, then $(\f_2)_*^{-1}\circ(\f_1)_*:\Gamma(E_1)\rightarrow\Gamma(E_2)$ is an isomorphism. A classical result given in  \cite[Proposition XIII p.78]{greub} states that the map $f\mapsto f_*$ is an isomorphism of $\con^\infty(M)$-modules from the set of morphisms from $E_1$ to $E_2$  to the set of morphisms from $\Gamma(E_1)$ to $\Gamma(E_2)$. Applying this result, there exists a unique isomorphism $h:E_1\rightarrow E_2$ such that $h_*=(\f_2)_*^{-1}\circ(\f_1)_*$. By construction, $(\f_2)_*\circ h_*=(\f_1)_*$ and applying again \cite[Proposition XIII p.78]{greub} we get $\f_2\circ h=\f_1$.
\hfill$\blacksquare$

Let $(E,\f)$ be a $(n,k)$-\rvd, $\bD=\f_*(\Gamma(E))=\{\f\circ\sigma\mid\sigma\in\Gamma(E)\}$ be its associated submodule and denote by $\bD(q)$  the linear subspace $\{V(q)\mid  V\in \bD\}=\f(E_q)\subseteq  T_q M$.
Let  $\mathrm{Lie}(\bD)$ be the smallest Lie subalgebra of
$\mathrm{Vec}(M)$
containing $\bD$ and, for every $q\in M$, let $\mathrm{Lie}_q(\bD)$ be the linear subspace of $T_qM$ whose elements are evaluation at $q$ of elements belonging to $\mathrm{Lie}(\bD)$.
We say that $(E,\f)$ satisfies the {\it Lie bracket generating condition} if
$\mathrm{Lie}_q(\bD)=T_q M$ for every $q\in M$.

A property $(P)$ defined for $(n,k)$-\rvd s is said to be {\it generic}
if for every vector bundle $E$ of rank $k$ over $M$, $(P)$ holds for every $\f$ in an  open and dense  subset of the set of  morphisms of vector bundles from $E$ to $TM$, endowed with the 
$\con^\infty$-Whitney topology. 
E.g., the Lie bracket generating condition is a generic property among  $(n,k)$-\rvd s satisfying $k>1$.

We say that a $(n,k)$-\rvd\ $(E,\f)$ is {\it orientable} if $E$ is orientable as a vector bundle.

A rank-varying sub-Riemannian structure is defined by requiring that $E$ is an Euclidean bundle.

\begin{definition}\label{gensrs}
A {\it $(n,k)$-\rvsR}  is a triple ${\cal S}=(E,\f,\langle\cdot,\cdot\rangle)$ where $(E,\f)$ is a Lie bracket generating $(n,k)$-\rvd\ on a  manifold $M$ and $\langle\cdot,\cdot\rangle_q$ is a   scalar product on $E_q$ smoothly depending on $q$. 
\end{definition}

Several classical structures can be seen as particular cases of \rvsR s, e.g.,  Riemannian structures and classical (constant-rank) sub-Riemannian structures (see \cite{book2,montgomery}). An $(n,n)$-\rvsR\ is called {\it $n$-dimensional almost-Riemannian structure}.  In this paper, we focus on  $2$-dimensional almost-Riemannian structures ($2$-ARSs for short). 
\bigskip

%%%%%%%%%%%%%%%%%%%%%%%%%%%%%%%%%%%%%%%%%%%%%%%%%%%%%%%%%%%%%%%%

Let ${\cal S}=(E,\f,\langle\cdot,\cdot\rangle)$ be a  $(n,k)$-\rvsR. The Euclidean structure on $E$ and the injectivity of the morphism $f_*$ allow to define a symmetric  positive definite  $\con^\infty(M)$-bilinear form on the submodule $\bD$ by
\begin{eqnarray*}
G:\bD\times\bD&\to&\con^\infty(M)\\
(V,W)&\mapsto&\langle\sigma_V,\sigma_W\rangle,
\end{eqnarray*}
where $\sigma_V,\sigma_W$ are the unique sections of $E$ such that $\f\circ\sigma_V=V, \f\circ\sigma_W=W$. 

If $\sigma_1,\dots,\sigma_k$ is an orthonormal frame for $\langle\cdot,\cdot\rangle$ on an open subset $\Omega$ of $M$, an {\it  orthonormal frame for $G$} on $\Omega$ is given by  $\f\circ\sigma_1,\dots,\f\circ\sigma_k$. Orthonormal frames are systems of local generators of $\bD$.

For every $q\in M$ 
and every $v\in\bD(q)$ define
\bqn
\Gq(v)=\inf\{\langle u, u\rangle_q \mid u\in E_q,\f(u)=v\}.
\eqnn

In this paper, a curve $\g:[0,T]\to M$  absolutely continuous with respect to the differential structure is said to be {\it admissible} for ${\cal S}$ 
if there exists a measurable essentially bounded function 
\bqn
[0,T]\ni t\mapsto u(t)\in E_{\g(t)}
\eqnn
called  {\it control function}, 
such that 
$\dot \g(t)=\f(u(t))$  for almost every $t\in[0,T]$. 
Given an admissible 
curve $\g:[0,T]\to M$, the {\it length of $\g$} is  
\bqn
\ell(\g)= \int_{0}^{T} \sqrt{ \gg_{\gamma(t)}(\dot \g(t))}~dt.\eqnn

The {\it Carnot-Caratheodory distance} (or sub-Riemannian distance) on $M$  associated with ${\cal S}$ is defined as
\bqn\nonumber
d(q_0,q_1)=\inf \{\ell(\g)\mid \g(0)=q_0,\g(T)=q_1, \g\ \mathrm{admissible}\}.
\eqn
The finiteness and the continuity of $d(\cdot,\cdot)$ with respect 
to the topology of $M$ are guaranteed by  the Lie bracket generating 
assumption on the \rvsR\ (see \cite{book2}).  
The Carnot-Caratheodory distance associated with ${\cal S}$  endows $M$ with the 
structure of metric space compatible with the topology of $M$ as differential manifold. 

We give now a characterization of admissible curves.

\begin{proposition}\label{lipadm}
Let $(E,f,\langle\cdot,\cdot\rangle)$ be a rank-varying sub-Riemannian structure on a manfold $M$.
Let $\gamma:[0,T]\rightarrow M$ be an absolutely continuous curve.
Then $\gamma$ is admissible if and only if it is Lipschitz continuous
with respect to
the sub-Riemannian distance.
\end{proposition}
{\bf Proof.} First we prove that if the curve is admissible  then
it is Lipschitz with respect to $d$ ({\it $d$-Lipschitz} for short). This is a direct consequence of the definition of the
sub-Riemannian distance. Indeed, let
\bqn
[0,T]\ni t\mapsto u(t)\in E_{\g(t)}
\eqnn
be a control function for $\gamma$ and  let $L>0$ be the  essential supremum of $\sqrt{\langle u,u\rangle}$.
Then,  for every subinterval $[t_0,t_1]\subset[0,T]$ one has
$$
d(\gamma(t_0),\gamma(t_1)) \leq \int_{t_0}^{t_1} \sqrt{ \gg_{\gamma(t)}(\dot \g(t))} dt \leq
\int_{t_0}^{t_1} \sqrt{\langle u(t),u(t)\rangle} dt \leq L(t_1-t_0).
$$
Hence $\gamma$ is $d$-Lipschitz.

Viceversa, assume that $\gamma$ is $d$-Lipschitz with Lipschitz constant $L$. Since $\gamma$ is absolutely continuous, it is differentiable almost everywhere on $[0,T]$.
Thanks to the Ball-Box Theorem (see \cite{bellaiche}), for every $t\in[0,T]$ such that the tangent vector $\dot \g(t)$ exists, $\dot\g(t)$ belongs to the distribution $\Delta(\gamma(t))$ (if not, the curve would fail to be $d$-Lipschitz). Hence for almost every $t\in[0,T]$ there exists  $u_t\in E_{\gamma(t)}$ such
that $\dot\gamma(t)=f(u_t)$.
Moreover, since the curve is $d$-Lipschitz,
one has that $\gg_{\g(t)}(\dot\gamma(t))\leq L^2$ for almost every $t\in [0,T]$. This can be seen computing lengths 
in
privileged coordinates (see \cite{bellaiche} for the definition of this system of coordinates). Hence, we can assume that $\langle u_t,u_t\rangle\leq L^2$ almost everywhere.
Finally, we apply Filippov Theorem (see \cite[Theorem 3.1.1 p.36]{bressan}) to the differential inclusion
$$
\dot\g(t)\in\{f(u)\mid \pi_E(u)=\gamma(t)\textrm{ and } \langle u,u\rangle \leq L^2\}.
$$
that assures the existence of a measurable choice of the control function corresponding to $\gamma$. Thus $\gamma$ is admissible.\hfill$\blacksquare$

Given a 2-ARS ${\mathcal
S}$, we define its {\it singular locus} as the set
$$\Zz=\{q\in M\mid \Delta(q) \subsetneq T_qM\}.$$ 
Since $\Delta$ is  bracket generating, the subspace $\bD(q)$ is nontrivial for every $q$ and $\Zz$ coincides with the set of points $q$ where $\bD(q)$ is one-dimensional.

We say that $\cal S$ {\it satisfies condition} \HH\  if the following properties hold:
{\bf (i)} $\Zz$ is an
embedded one-dimensional 
submanifold of
$M$;
{\bf (ii)} the points $q\in M$ at which $\bD_2(q)$ is
one-dimensional are isolated;
{\bf (iii)}  $\bD_3(q)=T_qM$ for every $q\in M$, where $\bD_1 = \bD$ and $\bD_{k+1}=\bD_k+[\bD,\bD_k]$.
It is not difficult to prove that property \HH\ is generic among 2-ARSs (see  \cite{ABS}). This hypothesis was essential to show Gauss--Bonnet type results for ARSs in \cite{euler,ABS,high-order}. The following  theorem recalls the  local normal forms for ARSs satisfying hypothesis \HH.
\begin{theorem}[\cite{ABS}]
\label{t-normal}
Given a  2-ARS ${\mathcal S}$ satisfiyng \HH, for every point
$q\in M$ there exist a neighborhood $U$ of $q$, an orthonormal frame
$(X,Y)$ for $G$ on $U$ and smooth coordinates defined on $U$ such that $q=(0,0)$ and $(X,Y)$
has one of the
forms
\bqn
\mathrm{(F1)}&& ~~X(x,y)=(1,0),~~~Y(x,y)=(0,e^{\phi(x,y)}), \nn  \\
\mathrm{(F2)}&& ~~X(x,y)=(1,0),~~~Y(x,y)=(0,x e^{\phi(x,y)}),\nn   \\
\mathrm{(F3)}&& ~~X(x,y)=(1,0),~~~Y(x,y)=(0,(y -x^2
\psi(x))e^{\pphi(x,y)}), \nn
\eqn
where $\phi$, $\pphi$ and $\psi$ are smooth real-valued functions such that
$\phi(0,y)=0$ and  $\psi(0)>0$.
\et
Let ${\mathcal S}$ be a 2-ARS satisfying \HH.
 A point $q\in M$ is said to be an
{\it ordinary point} if $\bD(q)=T_q M$, hence, if ${\mathcal
S}$ is locally described by (F1). We call $q$ a  {\it Grushin
point} if $\bD(q)$ is one-dimensional and $\bD_2(q)=T_q M$, i.e., if
the local description (F2) applies. Finally, if
$\Delta(q)=\Delta_2(q)$ has dimension one and $\bD_3(q)=T_q M$
then we say that $q$ is a {\it tangency point} and ${\mathcal
S}$ can be described near $q$ by the normal form (F3). We define
 $${\cal T}=\{q\in \Zz\mid q \mbox{ tangency point of } {\cal S}\}.$$
Assume ${\cal S}$ and $M$ to be oriented. Thanks to the hypothesis \HH, $M\setminus \Zz$  splits into two open sets
$M^+$ and $M^-$ such that $\f:E|_{M^+}\rightarrow TM^+$ is an orientation-preserving isomorphism and $\f:E|_{M^-}\rightarrow TM^-$ is an orientation-reversing isomorphism.

%%%%%%%%%%%%%%%%%%%%%%%%%%%%%%%%%%%%%%%%%%%%%%%%%%%%%%%%%%%%%%%%%

%%%%%%%%%%%%%%%%%%%%%%%%%%%%%%%%%%%%%%%%%%%%%%%%%%%%%%%%%%%
%%%%%%%%%%%%%%%%%%%%%%%%%%%%%%%%%%%%%%%%%%%%%%%%%%%%%%%%%%%%%

\section{Number of revolutions and graph of a $2$-ARS}\label{deftau}

From now on $M$ is a compact oriented surface and ${\cal S}=(E,\f,\langle\cdot,\cdot\rangle)$ is an oriented ARS 
on $M$ satisfying \HH.

 Fix on $\Zz$ the orientation induced by $M^+$ and consider a connected component $W$ of $\Zz$. Let $V\in\Gamma(TW)$ be a never-vanishing vector field  whose duality product with the fixed orientation on $W$ is positive.
Since $M$ is oriented, $TM|_{W}$ is isomorphic to the trivial bundle of rank $2$ over $W$.  We choose an isomorphism  $t:TM|_W\rightarrow W\times\R^2$ such that $t$ is orientation-preserving and  
for every $q\in W$, $t\circ V(q)=(q,(1,0))$. This trivialization induces an orientation-preserving isomorphism between the projectivization of $TM|_W$ and $W\times S^1$. 
For the sake of readability, in what follows we omit the isomorphism $t$ and identify $TM|_W$ (respectively, its projectivization) with $W\times \R^2$ (respectively, $W\times S^1$). 

 Since  $\bD|_W$ is a subbundle of rank one of $TM|_W$, $\bD|_W$ can be seen
 as a section of the projectivization of $TM|_W$, i.e., a smooth map (still denoted by $\bD$) $\bD:W\rightarrow W\times S^1$ such that $\pi_1\circ\bD=\mathrm{Id}_W$, where $\pi_1:W\times S^1\rightarrow W$ denotes the projection on the first component. 
We define $\tau(\bD,W)$, the {\it number of \rotations\ } of $\bD$ along $W$, to be the degree of the map $\pi_2\circ\bD:W\rightarrow S^1$, where $\pi_2:W\times S^1\rightarrow S^1$ is the projection on the second component. Notice that $\tau(\bD,W)$ changes sign if we reverse the orientation of $W$.

Let us show how to compute $\tau(\bD,W)$. By construction, $\pi_2\circ  V:W\rightarrow S^1$ is constant. 
Let $\pi_2\circ V(q)\equiv \th_0$. Since $\bD_3(q)=T_qM$ for every $q\in M$,  $\th_0$ is  a regular value of $\pi_2\circ \bD$. 
 By definition,
\begin{equation}\label{mario}
\tau(\bD,W)=\sum_{q\mid\pi_2\circ\bD(q)=\th_0}\sign( d_q(\pi_2\circ\bD))=\sum_{q\in W\cap {\cal T}}\sign(d_q(\pi_2\circ\bD)),
\end{equation}
where $d_q$ denotes the differential at $q$ of a smooth map and $\sign( d_q(\pi_2\circ\bD))=1,$ resp. $-1$, if $d_q(\pi_2\circ\bD)$ preserves, resp. reverses, the orientation. The equality in \r{mario} follows from the fact that a  point $q$ satisfies $\pi_2\circ\bD(q)=\th_0$ if and only if $\bD(q)$ is tangent to $W$ at $q$, i.e., $q\in {\cal T}$.

  Define the {\it contribution at a tangency point $q$} as $\tau_q=\sign(d_q(\pi_2\circ\bD))$ (see Figure~\ref{tauu}). Moreover,  we define
$$
\tau({\cal S})=\sum_{W\in\CC(\Zz)}\tau(\bD,W),
$$ 
where $\CC(\Zz)=\{W\mid W \textrm{ connected component of } \Zz\}$. Clearly, $\tau({\cal S})=\sum_{q\in{\cal T}}\tau_q$.

\begin{figure}[h!]
\begin{center}
\input{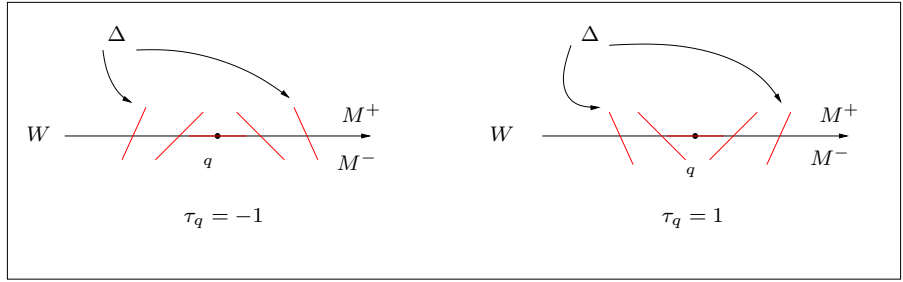}

\caption{Tangency points with opposite contributions}\label{tauu}
\end{center}
\end{figure}

\bigskip

Let us associated with the 2-ARS ${\cal S}$
the graph ${\cal G}=({\cal V}(\Gg),{\cal E}(\Gg))$ where
\bi
\iii each vertex in ${\cal V}(\Gg)$ represents a connected component of $M\setminus \Zz$;
\iii each edge in ${\cal E}(\Gg)$ represents a connected component of $\Zz$; 
\iii the edge corresponding to a connected component $W$ connects the two vertices corresponding to the connected components $M_1$ and $M_2$ of $M\setminus \Zz$ such that $W\subset \partial{M}_1 \cap \partial{M}_2$.
\ei

Thanks to the hypothesis \HH, every connected component of $\Zz$ joints a connected component of $M^+$ and one of $M^-$. Thus the graph ${\cal G}$ turns out to be {\it bipartite}, i.e., there exists a partition of the set of vertices into two subsets $V^+$ and $V^-$
such that each edge of ${\cal G}$ joins a vertex of $V^+$ to a vertex of $V^-$. 
Conversely, it is not difficult to see that every finite bipartite graph can be obtained from an oriented 2-ARS (satisfying \HH) on a compact oriented surface.

Using the bipartite nature of $\Gg$ we introduce an orientation on $\Gg$ given by  two functions $\alpha,\omega:{\cal E}(\Gg)\to {\cal V}(\Gg)$ defined as follows. If $e$ corresponds to $W$ then $\alpha(e)=v$ and  $\omega(e)=w$,  where $v$ and $w$ correspond respectively to the connected components $M_v\subset M^-$  and $M_w\subset M^+$ such that $W\subseteq \partial M_v\cap \partial M_w$.

We label each vertex $v$ corresponding to a connected component $\hat M$ of $M\setminus \Zz$ 
with a pair $(\sign(v),\chi(v))$ where sign$(v)=\pm1$ if $\hat M\subset M^\pm$ 
 and $\chi(v)$ is the Euler characteristic of $\hat M$. 
We define for every $e\in E({\cal G})$ the number $\tau(e)=\sum_{q\in W\cap {\cal T}}\tau_q$, where $W$ is the connected component of $\Zz$ corresponding to $e$.

Finally, we define a label  for each edge $e$ corresponding to a connected component $W$ of $ \Zz$ containing tangency points. 
Let $s\geq1$ be the cardinality of the set $W\cap {\cal T}$. The label of $e$ is an equivalence class of $s$-uples with entries in $\{\pm 1\}$ defined as follows. Fix on $W$ the orientation induced by $M^+$ and choose a point $q\in W\cap{\cal T}$. Let $q_1=q$ and for every $i=1,\dots, s-1$ let $q_{i+1}$ be the first element in $W\cap {\cal T}$ that we meet after $q_i$ walking  along $W$ in the fixed orientation. 
We associate with $e$ the equivalence class of $(\tau_{q_1},\tau_{q_2},\dots,\tau_{q_s})$ in the set of $s$-uples with entries in $\{\pm 1\}$ modulo cyclic permutations. In figure \ref{figsez3} an ARS on a surface of genus 4 and its labelled graph (figure \ref{figsez3}(a)) are portrayed. According to our definition of labels on edges, figures  \ref{figsez3}(a) and \ref{figsez3}(b) represent equal graphs associated with the same ARS. On the other hand,  the graph in figure \ref{figsez3}(c) is not the graph associated to the ARS of  figure \ref{figsez3}.  In figure \ref{come} two steps in the construction of the labelled graph associated with the ARS in figure \ref{figintro} are shown.
 
\brem Once an  orientation on $E$ is fixed  the labelled graph associated with ${\cal S}$ is unique. 
\erem
\begin{figure}[h!]
$$\input{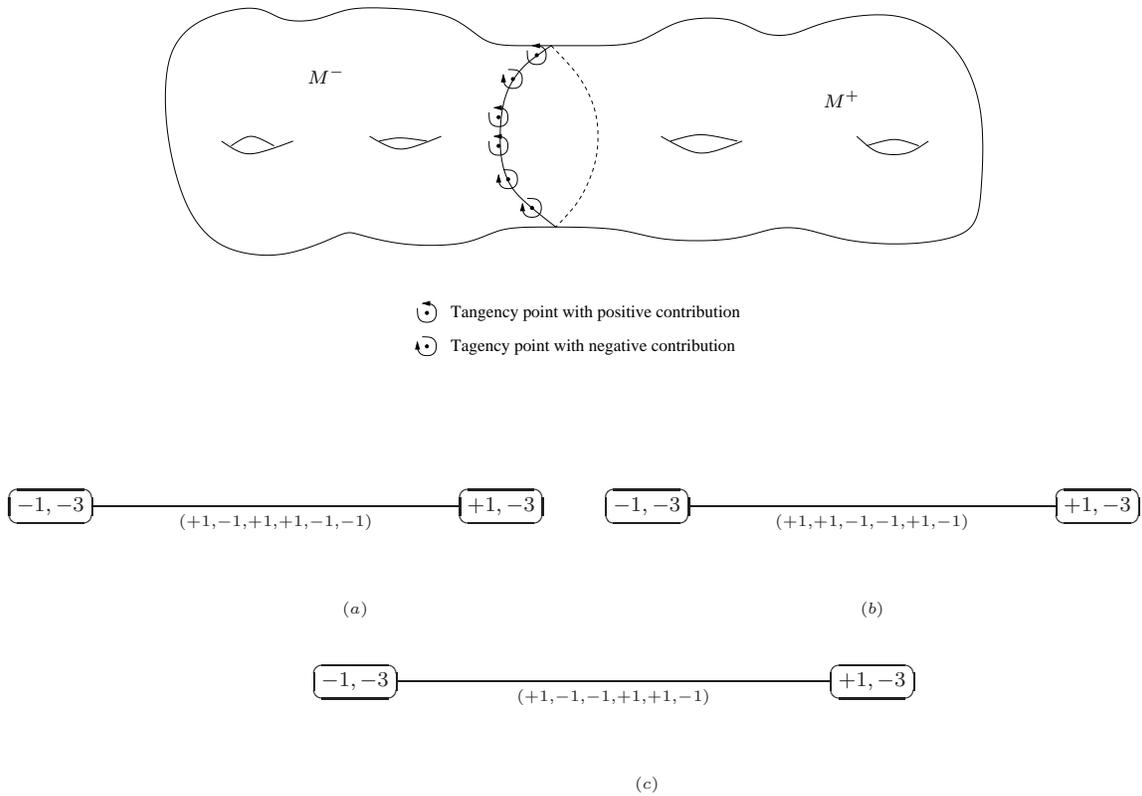}$$

\scriptsize{
$$
\xymatrix{&&\\
\ver{-1,-3}\ar@{-}[rrrr]_{(+1,-1,+1,+1,-1,-1)}&&&&\ver{+1,-3}&
\ver{-1,-3}\ar@{-}[rrrr]_{(+1,+1,-1,-1,+1,-1)}&&&&\ver{+1,-3}
\\
&&&\ar@{}_{(a)}&&&&\ar@{}_{(b)}\\
&&&\ver{-1,-3}\ar@{-}[rrrr]_{(+1,-1,-1,+1,+1,-1)}&&&&\ver{+1,-3}\\
&&&&&\ar@{}_{(c)}
}
$$
}
\caption{\label{figsez3} Example of ARS  on a surface of genus 4. Figures (a) and (b) illustrate equal labelled graphs associated with the ARS. Figure (c) gives an example of labelled graph different from the graph in figure (a)}
\end{figure}

\begin{figure}
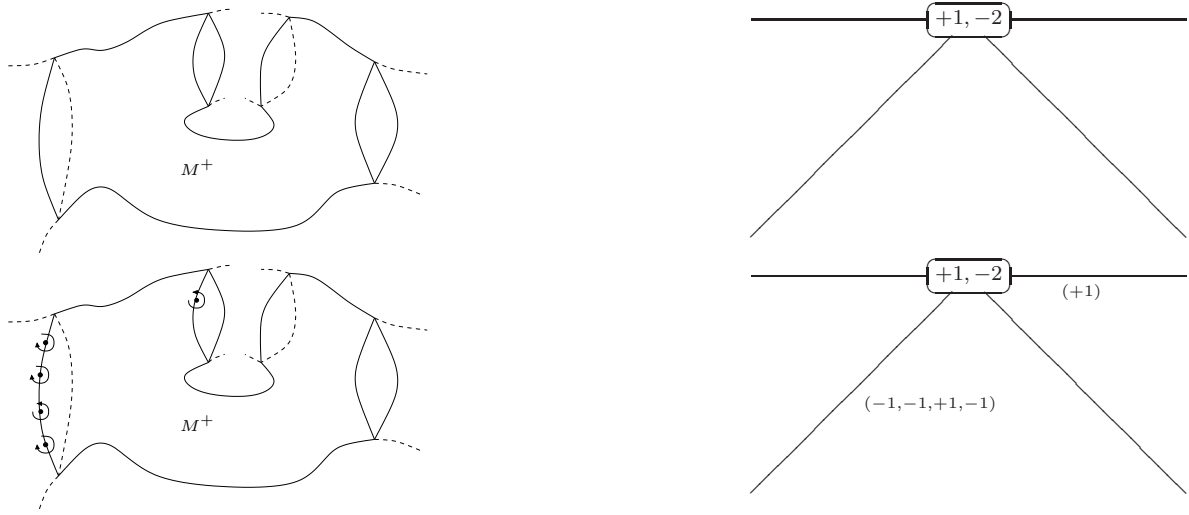

\scriptsize{
\xymatrix@=3cm@!0{*!<4cm,1.5cm>{\input{punto1.pstex_t}}&&\ver{+1,-2}\ar@{-}[l]\ar@{-}[r]\ar@{-}[dr]\ar@{-}[dl]&\\
&&&\\
}
\xymatrix@=3cm@!0{*!<4cm,1.5cm>{\input{punto2.pstex_t}}&&\ver{+1,-2}\ar@{-}[l]\ar@{-}[r]_{(+1)}\ar@{-}[dr]\ar@{-}[dl]^{(-1,-1,+1,-1)}&\\
&&&}
}
\caption{\label{come}Algorythm to build the graph}
\end{figure}

We define  an equivalence 
relation on the set of graphs associated with oriented ARS on $M$ satisfying hypothesis \HH.

\bdeff\label{samegraph}
Let ${\cal S}_i=(E_i,\f_i, \langle\cdot,\cdot\rangle_i)$, $i=1,2$,
 be two oriented almost-Riemannian structures on a compact
 oriented surface $M$ satisfying hypothesis \HH. Let ${\cal G}_i$ be the labelled graph associated with ${\cal S}_i$ and denote by
$\alpha_i,\omega_i:\Ee(\Gg_i)\rightarrow\V(\Gg_i)$ the functions defined as above. We say that 
${\cal S}_1$ and ${\cal S}_2$ have  \emph{\same\ graphs} if , after possibly changing the orientation on $E_2$, they have the same labelled graph. 

In other words,  after possibly changing the orientation on $E_2$ and still denoting by ${\cal G}_2$ the associated graph, there exist bijections $u:\V({\cal G}_1)\rightarrow \V({\cal G}_2)$, $k:\Ee({\cal G}_1)\rightarrow \Ee({\cal G}_2)$ such that  
 the diagram 
\begin{equation}\label{d1}
\xymatrix{
 \V(\Gg_1)  \ar[r]_u   & \V(\Gg_2)             \\
\Ee(\Gg_1) \ar[u]^{\alpha_1}\ar[r]^k      &\Ee(\Gg_2)\ar[u]_{\alpha_2}                      
}    
\end{equation}
commutes and $u$ and $k$ preserve labels.
\edeff
Figure  \ref{grafoequivalente} illustrates the graph associated with the ARS obtained by reversing the orientation of the ARS in figure \ref{figintro}.

\begin{figure}[h!]
\scriptsize{
$$\xymatrix{
&&\ver{+1,-2}\ar@{-}[ddll]\ar@{-}[ddrr]_{(+1,-1,+1,+1)}&&&\ver{+1, 0}\ar@{-}@(ul,dl)[ddl]_{(-1)}\ar@{-}[ddl]&&\ver{+1,-4}\ar@{-}[ddlll]\ar@{-}[ddl]^{(+1,+1)}\ar@{-}[ddr]\ar@{-}@(ur,dr)[ddr]&\\
&&&&&&\\
\ver{-1,+1}&&&&\ver{-1,-2}&&\ver{-1,+1}&&\ver{-1, 0}
}
$$}

\caption{\label{grafoequivalente}Equivalent graph to the one in figure \ref{figintro}}
\end{figure}
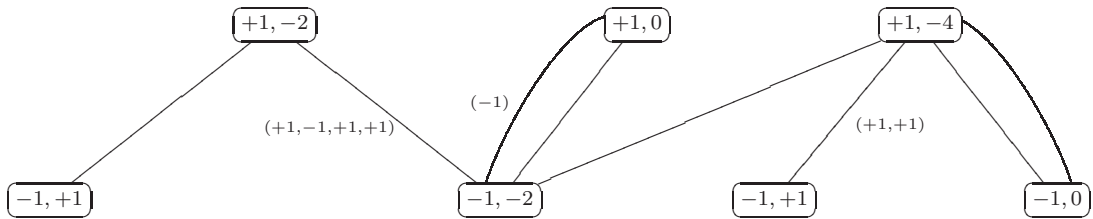

\section{Lipschitz equivalence}\label{proof}

This section is devoted to the proof of Theorem~\ref{lip-eq} which is a generalization to ARSs of the well-known fact that all Riemannian structures on a compact oriented surface are Lipschitz equivalent.

Let $M_1, M_2$ be two manifolds. For $i=1,2$, let  ${\cal S}_i=(E_i,\f_i,\langle\cdot,\cdot\rangle_i)$ be a sub-Riemannian structure on $M_i$. Denote by
$d_i$ the Carnot--Caratheodory distance on $M_i$ associated with ${\cal S}_i$.
\bdeff\label{deflipeq} 
We say that a diffeomorphism $\PH:M_1\rightarrow M_2$ is a \emph{Lipschitz equivalence} 
if it is bi-Lipschitz as a map from $(M_1,d_1)$ to $(M_2,d_2)$.
\edeff

Notice that in Theorem~\ref{lip-eq} we can assume  $M_1=M_2=M$. Indeed, if two ARSs are Lipschitz equivalent, then by definition there exists a diffeomorphism $\varphi:M_1\rightarrow M_2$. On the other hand, if the associated graphs are equivalent, by \cite[Theorem 1]{euler} it follows that $E_1$ and $E_2$ are isomorphic vector bundles. Hence the underlying surfaces are diffeomorphic.

\subsection{Necessity}\label{necessity}

Denote by $M_{i}^+$, respectively $M_i^-$, the set where $\f_i$ is an orientation-preserving,
 respectively orienta\-tion-reversing, isomorphism of vector bundles, and by $\bD^i$ the submodule $\{\f_i\circ\sigma\mid\sigma\in\Gamma(E_i)\}$. Let $\Zz_i$ be the singular locus 
of ${\cal S}_i$ and ${\cal T}_i$ the set of tangency points of ${\cal S}_i$. Finally,  for every
 $q\in {\cal T}_i$, denote by $\tau^i_q$ the contribution at the tangency point defined in Section~\ref{deftau} with $\bD=\bD^i$.
     
In this section we  assume $\PH:(M,d_1)\rightarrow (M,d_2)$ to be a Lipschitz equivalence and we show that ${\cal S}_1$ and
 ${\cal S}_2$ have  \same\ graphs. As a consequence of the Ball-Box Theorem (see, for instance, \cite{bellaiche})
 one can prove the following 
result. 

\bl\label{samepoint}
If $p$ is an ordinary, Grushin or tangency
 point for ${\cal S}_1$, then $\PH(p)$ is an ordinary, Grushin or tangency
 point for ${\cal S}_2$, respectively.     
\el

Thanks to Lemma~\ref{samepoint}, for every connected component $\hat M$  of $M\setminus \Zz_1$,  
$\PH(\hat M)$ is a connected component of $M\setminus \Zz_2$  and for every connected component $W$ of 
$\Zz_1\cap\partial\hat M$, $\PH(W)$ is a connected component of $\Zz_2\cap \partial\PH(\hat M)$. 
 Moreover, since $\PH|_{\overline{\hat M}}$ 
is a diffeomorphism, it follows that $\chi(\hat M)=\chi(\PH(\hat M))$. After possibly changing the orientation
on $E_2$, we may assume $\PH(M_1^\pm)=M_{2}^\pm$. We will prove that, in this case, the labelled graphs are equal. Indeed, if $v\in \V({\cal G}_1)$ corresponds to $\hat M$, define $u(v)\in \V({\cal G}_2)$ as the vertex 
corresponding to $\PH(\hat M)$. If $e\in \Ee({\cal G}_1)$ corresponds to $W$ define $k(e)\in \Ee({\cal G}_2)$ as the edge 
corresponding to $\PH(W)$.
Then $\chi(u(v))=\chi(v)$, $\sign(u(v))=\sign(v)$, and, by construction, the diagram \r{d1} commutes.

 Let us compute the contribution at a tangency point $q$ of an ARS $(E,f,\langle\cdot,\cdot\rangle)$ using the corresponding normal form given in Theorem \ref{t-normal}.
\bl\label{calcotau}
 Let $\gamma:[0,T]\rightarrow M$ be a smooth 
 curve such that $\gamma(0)=q\in {\cal T}$ and $\dot\gamma(0)\in\bD(q)\setminus\{ 0\}$. Assume moreover that $\gamma$ is $d$-Lipschitz, where $d$ is the almost-Riemannian distance, and that $\gamma((0,T))$ is contained in one of the two connected components of $M\setminus\Zz$. Let $(x,y)$ be a  coordinate system centered at $q$ such that the form (F3) of Theorem \ref{t-normal} applies. Then $\gamma((0,T))\subset \{(x,y)\mid y-x^2\psi(x)<0\}$. Moreover, if  $\{(x,y)\mid y-x^2\psi(x)<0\}\subseteq M^+$, resp. $M^-$, then $\tau_q=1$, resp. $-1$. 
\el
{\bf Proof.} Since $\g(0)=(0,0)$ and $\dot\g(0)\in\mathrm{span}\{(1,0)\}\setminus\{0\}$, 
 there exist two smooth  functions $\ov{ x}(t),\ov{y}(t)$ such that $\gamma(t)=(t\ov{x}(t),t^2\ov{y}(t))$ and $\ov{x}(0)\neq 0$. Assume by contradiction that $\gamma((0,T))\subset \{(x,y)\mid y-x^2\psi(x)>0\}$, i.e., for $t\in(0,T)$, $\ov{y}(t)>\psi(t \ov{x}(t))\ov{x}(t)^2$. Since $\psi(0)>0$, for $t$ sufficiently small $\psi(t\ov{x}(t))>0$ and   $\ov{y}(t)^{1/3}>\psi(t\ov{x}(t))^{1/3}|\ov{x}(t)|^{2/3}$.   By the Ball-Box Theorem (see \cite{bellaiche}) there exist $c_1,c_2$ positive constants such that, for $t$ sufficiently small  we have
$$
c_1(|t \ov{x}(t)|+|t^2\ov{y}(t)|^{1/3})\leq d(\gamma(t),(0,0))\leq c_2(|t \ov{x}(t)|+|t^2\ov{y}(t)|^{1/3}).
$$
On the other hand, for $t$ sufficiently small,
$$
|t \ov{x}(t)|+|t^2\ov{y}(t)|^{1/3}>t^{2/3}|\ov{x}(t)|^{2/3}\psi(t \ov{x}(t))^{1/3}.
$$
Hence, for $t$ sufficiently small, $d(\gamma(t),(0,0))>c_3t^{2/3}$, with $c_3>0$.
This implies that $\gamma$ is not Lipschitz with respect to the almost-Riemannian distance. 
Finally, a direct computation shows the assertion concerning $\tau_q$, see Figure \ref{tauu}.
\hfill$\blacksquare$

 Next lemma, jointly with Lemma~\ref{samepoint}, guarantees that the two bijections $u$ and $k$ preserve labels.
 \bl\label{prod}
Let $q\in{\cal T}_1$. Then $\tau^{1}_q=\tau^{2}_{\PH(q)}$. 
\el
{\bf Proof.}
Apply Theorem~\ref{t-normal} to  ${\cal S}_1$ and find a neighborhood $U$  of $q$ and a coordinate system $(x,y)$ on $U$  such that $q=(0,0)$ and 
$\Zz_1\cap U=\{(x,y)\mid y=x^2\psi(x)\}$. Let $\sigma,\rho\in\Gamma(E|_U)$ be the local orthonormal frame such that $\f_1\circ\sigma=X$ and $f_1\circ\rho=Y$. Assume that $U_1^+=M_1^+\cap U=\{(x,y)\mid y-x^2\psi(x) >0\}$. Fix $T>0$ and consider the smooth curve $\gamma:[0,T]\rightarrow U$ defined by $\g(t)=(t,0)$. Then $\g$ is  admissible  for ${\cal S}_1$ with control function $u(t)=\sigma(t,0)$. By definition, for $T$ sufficiently small $\gamma((0, T))$ lies in a single connected component of $U\setminus \Zz_1$. Moreover, by Proposition~\ref{lipadm},  $\gamma$ is  a $d_1$-Lipschitz map with Lipschitz constant less or equal to $1$. Hence, according to Lemma~\ref{calcotau}, $\tau^1_q=-1$. 

Consider the curve $\tilde\gamma=\varphi\circ\gamma:[0,T]\rightarrow \varphi(U)$. Since $\varphi$ is Lipschitz, $\tilde\gamma$ is $d_2$-Lipschitz as a map from the  interval $[0,T]$ to the metric space $(\varphi(U),d_2)$.  Moreover, $\tilde\gamma$ is smooth and $\dot{\tilde{\gamma}}(0)\in\bD^2(\varphi(q))\setminus\{ 0\}$, $\varphi$ being a diffeomorphism mapping $\Zz_1$ to $\Zz_2$. Finally, since $\varphi(M_1^-)=M_2^-$, then  $\tilde{\gamma}((0,T))\subset U_2^-=\varphi(U)\cap M_2^-$.  
 Thus, by Lemma~\ref{calcotau}, $\tau^2_{\varphi(q)}=-1$. Analogously, one can prove the statement in the case $U_1^+=\{(x,y)\mid y-x^2\psi(x) <0\}$ (for which $\tau_q^1=\tau_q^2=1$).
\hfill$\blacksquare$

 Lemma~\ref{prod} implies  that 
${\cal S}_{1}$ and ${\cal S}_2$ have equal labelled graphs.
This concludes the proof that having  \same\ graphs is a necessary condition for two ARSs being Lipschitz equivalent.

\subsection{Sufficiency}\label{sufficiency}
In this section we prove that if ${\cal S}_{1}$ and ${\cal S}_2$ have \same\ graphs then there exists a Lipschitz 
equivalence between $(M,d_{1})$ and $(M,d_2)$.
After possibly changing the orientation on $E_2$, we assume the associated labelled graphs to be equal, i.e.,  there exist two bijections $u$, $k$ as in Definition~\ref{samegraph} such that diagram \r{d1} commutes.  

The proof is in five steps. 
The first step consists in proving  that we may assume $E_1=E_2$.
The second step  shows that we can restrict to the case $\Zz_1=\Zz_2$ and ${\cal T}_1={\cal T}_2$.
In the third step we prove that we can assume that $\bD^1(q)=\bD^2(q)$ at each point $q\in M$.
As fourth step, we demonstrate that the submodules $ \bD^1$ and $\bD^2$ coincide. In the fifth and final step we remark that we can assume $f_1=f_2$ and conclude.
The Lipschitz equivalence between the two structures will be the composition of the  diffeomorphisms singled out  in steps 1, 2, 3, 5.

 By construction, the push-forward of ${\cal S}_1$ along a diffeomorphism  $\psi$ of $M$, denoted by $\psi_*{\cal S}_1$, is Lipschitz equivalent to ${\cal S}_1$
and   has the
same labelled graph of ${\cal S}_{1}$.  
Notice, moreover, that the singular locus  of $\psi_*{\cal S}_1$ coincides with $\psi(\Zz_1)$
 and the set of tangency points coincides with $\psi({\cal T}_1)$.

{\bf Step 1.}
Having the same labelled graph implies 
$$
\sum_{v\in \V({\cal G}_{1})}\sign(v)\chi(v)+\sum_{e\in \Ee({\cal G}_{1})}\tau(e)=\sum_{v\in \V({\cal G}_2)}\sign(v)\chi(v)+\sum_{e\in \Ee({\cal G}_2)}\tau(e).
$$
By \cite[Theorem 1]{euler}, this is equivalent to say that the Euler numbers of $E_{1}$ and $E_2$ are equal. Since $E_{1}$ and $E_2$ are
oriented vector bundles of rank $2$,
 with the same Euler number, over a compact oriented surface, then they are isomorphic. Hence, we assume 
$E_1=E_2= E$.

{\bf Step 2.}
Using the bijections $u,k$ and the classification of compact oriented 
surfaces with boundary (see, for instance, \cite{hirsch}), one can prove the following lemma. 

\bl
There exists a diffeomorphism 
$\tilde{\PH}:M\rightarrow M$ such that  $\tilde\PH(M_{1}^+)=M_2^+$, $\tilde\PH(M_{1}^-)=M_2^-$, $\tilde\PH|_{\Zz_{1}}:\Zz_{1}\rightarrow \Zz_2$ is a diffeomorphism 
that maps ${\cal T}_{1}$ 
into ${\cal T}_2$,  and, for every $q\in{\cal T}_{1}$ $\tau^2_{\tilde\PH(q)}=\tau^1_q$.
 Moreover, if $v\in\V(\Gg_1)$ corresponds to  $\hat M\subset M\setminus \Zz_1$, then $\tilde\PH(\hat M)$ is the connected component of $M\setminus \Zz_2$ corresponding to $u(v)\in\V(\Gg_2)$;  if $e\in\Ee(\Gg_1)$ corresponds to
$W\subset \Zz_1$, then $\tilde\PH(W)$ is the connected component of $\Zz_2$ corresponding to $k(e)\in\Ee(\Gg_2)$.
 \el

The lemma implies that the singular locus  of $\tilde{\PH}_*{\cal S}_1$ coincides with $\Zz_2$
 and the set of tangency points coincides with ${\cal T}_2$. For the sake of readability, in the following we rename $\tilde{\PH}_*{\cal S}_1$ simply by ${\cal S}_1$ and we will denote by $\Zz$ the singular locus of the two structures, by ${\cal T}$ the set of their tangency
points, and by $M^\pm$ the set $M_i^\pm$. 
 
{\bf Step 3.}
Remark that the subspaces $ \bD^{1}(q)$ and $\bD^2(q)$ 
coincide at every ordinary and tangency point $q$. We are going to show that there exists a diffeomorphism
of $M$ that carries $\bD^1(q)$ into  $\bD^2(q)$ at every point $q$ of the manifold.
 
\bl\label{matrice}
Let $W$ be a  connected component of $\Zz$. There exist a tubular neighborhood $\cil$
 of $W$ and a diffeomorphism
$\PH_W:\cil\rightarrow \PH_W(\cil)$ such that $d_q\PH_W(\bD^1(q))=\bD^2(\PH_W(q))$ for every $q\in\cil$,  
$\PH_W|_W=\Id|_W$ and $\PH(\cil\cap M^\pm)\subset M^\pm$.
\el

{\bf Proof.} 
The idea of the proof is first to consider a smooth section  $A$
of Hom$(TM|_W;TM|_W)$ such that for every $q\in W$, $A_q:T_qM\rightarrow T_qM$ is an isomorphism and $A_q(\bD^1(q))=\bD^2(q)$. Secondly, we build a diffeomorphism $\varphi_W$ of a tubular neighborhood of $W$ such that $d_q\varphi_W= A_q$ for every point $q\in W$.

Choose on a tubular neighborhood $\cil$ of $W$ 
a parameterization $(\th,t)$ such that $W=\{(\th,t)\mid t=0\}$, $M^+\cap \cil=\{(\th,t)\mid t>0\}$ and
$\frac{\partial}{\partial \th}\left|_{(\th,0)}\right.$ induces on $W$ the same orientation as $M^+$. We are going to show the existence of two smooth functions $a,b:W\rightarrow \R$ such that $b$ is positive and for every $(\th,0)\in W$, 
\begin{equation}\label{endom}
\left(\ba{cc}1 & a(\th)\\0&b(\th)\ea\right)(\bD^1(\th,0))=\bD^2(\th,0).
\end{equation}
Then, for every $q=(\th,0)\in W$ defining $A_q:T_qM\rightarrow T_qM$ by
\begin{equation}\label{defisom}
A_{(\th,0)}=\left(\ba{cc}1 & a(\th)\\0&b(\th)\ea\right),
\end{equation}
we will get an isomorphism smoothly depending on the point $q$ and carrying $\bD^{1}(q)$ into $\bD^{2}(q)$.

 Let $W\cap{\cal T}=\{(\th_1,0),\dots,(\th_s,0)\}$, with $s\geq 0$.
Using the chosen parametrization, there exist two smooth functions $\beta_1,\beta_2:W\setminus\{(\th_1,0),\dots,(\th_s,0)\}\rightarrow \R$ such that
$\bD^i(\th,0)=\mathrm{span}\{(\beta_i(\th),1)\}$. For every $j=1,\dots s$, there exists a smooth function $g^i_j$ defined on a neighborhood of $(\th_j,0)$ in $W$ such that $g^i_j(\th_j)\neq 0$, $\tau^i_{(\th_j,0)}=\sign(g^i_j(\th_j))$ and
$$
\beta_i(\th)=\frac{1}{(\th-\th_j)g^i_j(\th)}, \quad\th\sim\th_j.
$$
Since the graphs associated with ${\cal S}_1, {\cal S}_2$ are equivalent, for every $j=1\dots s$  we have $\tau^1_{(\th_j,0)}=\tau^2_{(\th_j,0)}$. Hence $\frac{g^2_j(\th_j)}{g^1_j(\th_j)}>0$ for every $j$.
Let $b:W\rightarrow \R$ be a positive smooth function such that for each $j\in\{1,\dots s\},\,b(\th_j)=\frac{g^2_j(\th_j)}{g^1_j(\th_j)}$. Define $a:W\rightarrow\R$ by
$$
a(\th)=b(\th)\beta_2(\th)-\beta_1(\th).
$$
Clearly $a$ is smooth on $W\setminus\{(\th_1,0),\dots,(\th_s,0)\}$. Moreover, thanks to our choice of $b$,
$a$ is smooth at $\th_j$, and, by construction, we have \r{endom}. The existence of $a,b$ is established.

Define $A_q$ as in \r{defisom}.
Let us extend the isomorphism $A_q$ defined for $q\in W$ to a tubular neighborhood. Define $\PH_W:\cil\rightarrow \cil$  by
$$
\PH_W(\th,t)=(a(\th)t+\th,b(\th)t).
$$ 
By construction,  $d_{(\th,0)}\PH_W$ is an isomorphism. Hence, reducing $\cil$ if necessary, 
$\PH_W:\cil\rightarrow \PH_W(\cil)$
turns out to be a diffeomorphism. Finally, by definition, $\PH_W(\th,0)=(\th,0)$ and, since $b$ is positive, $\PH(\cil\cap M^\pm)\subset M^\pm$.
\hfill$\blacksquare$

We apply Lemma~\ref{matrice} to every connected component $W$ of $\Zz$. We reduce, 
if necessary,  the tubular neighborhood $\cil$ of $W$ in such a way that every pair of distinct connected 
component of $\Zz$ have disjoint corresponding tubular neighborhoods built as in Lemma~\ref{matrice}. 
We claim that
there exists a diffeomorphism $\PH:M\rightarrow M$ such that $\PH|_\cil=\PH_W$ for every connected component $W$ 
of $\Zz$. This is a direct consequence of the fact that the labels on vertices of ${\cal G}_1$ and ${\cal G}_2$ 
are equal and of the classification of compact oriented surfaces with boundary (see \cite{hirsch}). 
By construction, the 
push-forward of ${\cal S}_1$ along $\PH$  is Lipschitz equivalent to ${\cal S}_1$ and has the same labelled graph as
${\cal S}_1$. To simplify notations, we denote $\PH_*{\cal S}_1$ by ${\cal S}_1$.   By Lemma~\ref{matrice}, $\bD^1(q)=\bD^2(q)$
at every point $q$.  

{\bf Step 4.}
The next point is to prove that
 $ \bD^1$ and $\bD^2$ coincide as $\con^\infty(M)$-submodules. 

\bl\label{coin}
The submodules $\bD^1$ and $\bD^2$ associated with ${\cal S}_1$ and ${\cal S}_2$ coincide.
\el

{\bf Proof.}
It is sufficient to show that
 for every $p\in M$ there exist a neighborhood $U$ of $p$ such that $\bD^1|_U$ and 
$\bD^2|_U$ are generated as $\con^\infty(M)$-submodules by the same pair of vector fields.

If $p$ is an ordinary point, then taking  $U=M\setminus \Zz$, we have $\bD^1|_U=\bD^2|_U=\mathrm{Vec}(U)$.

Let $p$ be a Grushin point  and apply Theorem~\ref{t-normal} to ${\cal S}_1$ to find a neighborhood  $U$  of $p$  such that
$$
\bD^1|_U=\mathrm{span}_{\con^\infty(M)}\{F_1,F_2\}, ~\mathrm{ where }~~ F_1(x,y)=(1,0),~F_2(x,y)=(0,x e^{\phi(x,y)}).
$$
Up to reducing $U$ we assume the existence of a frame 
$$
G_1(x,y)=(a_1(x,y),a_2(x,y)),~~~G_2(x,y)=(b_1(x,y),b_2(x,y))
$$
 such that $\bD^2|_U=\mathrm{span}_{\con^\infty(M)}\{G_1,G_2\}$. Since $\bD^1(q)=\bD^2(q)$ at every point $q\in M$,  $a_2(0,y)\equiv 0$ and $b_2(0,y)\equiv 0$. Since $\bD^2(0,y)$ is one-dimensional, let us assume $a_1(0,y)\neq 0$ for every $y$. Moreover,  after possibly further reducing $U$, $\bD^2|_U=\mathrm{span}_{\con^\infty(M)}\{(1/a_1)G_1, G_2-(b_1/a_1)G_1\}$ hence we  may assume 
$a_1(x,y)\equiv 1$ and $b_1(x,y)\equiv 0$.  
The conditions $a_2(0,y)\equiv 0$ and $b_2(0,y)\equiv 0$  imply $a_2(x,y)=x \overline{a}_2(x,y)$ and  $b_2(x,y)=x \overline{b}_2(x,y)$ respectively, with $\overline{a}_2,\overline{b}_2$ smooth functions.
Since $[G_1,G_2]|_{(0,y)}=(0,\overline{b}_2(0,y))$, thanks to 
hypothesis \HH\ on ${\cal S}_2$, we have $\overline{b}_2(0,y)\neq 0$. Hence, reducing $U$ if necessary,
\begin{eqnarray}
\bD^2|_U&=&\mathrm{span}_{\con^\infty(M)}\{G_1-(\overline{a}_2(x,y)/\overline{b}_2(x,y))G_2, (e^{\phi(x,y)}/\overline{b}_2(x,y)) G_2\}\nonumber\\
&=&\mathrm{span}_{\con^\infty(M)}\{F_1,F_2\}=\bD^1|_U.\nonumber
\end{eqnarray}

Finally, let $p$ be a tangency point. Apply Theorem~\ref{t-normal} to ${\cal S}_1$, i.e., choose  a neighborhood $U$ of $p$
and a system of coordinates $(x,y)$ such that $p=(0,0)$,
$$
\bD^1|_U=\mathrm{span}_{\con^\infty(M)}\{F_1,F_2\}, ~\mathrm{ where }~~ F_1(x,y)=(1,0),~F_2(x,y)=(0,(y-x^2\psi(x))e^{\pphi(x,y)}),
$$ 
and $\psi,\pphi$ are smooth functions such that $\psi(0)>0$.
Consider the change of coordinates
$$
\tilde x=x,~~~\tilde y=y-x^2\psi(x).
$$
Then 
$$
F_1(\tilde x,\tilde y)=(1,\tilde x a(\tilde x)),~~~F_2(\tilde x,\tilde y)=(0, \tilde y e^{\pphi(\tilde x,\tilde y+\tilde x^2\psi(\tilde x))}),
$$
where $a(\tilde x)=-2\psi(\tilde x)-\tilde x\psi'(\tilde x)$.
To simplify notations, in the following we rename $\tilde x, \tilde y$ by $x,y$ respectively and we still denote by $\pphi(x,y)$ the function $\pphi(x,y+x^2\psi(x))$.  In the new coordinate system we have  $p=(0,0)$, $\Zz\cap U=\{(x,y)\mid y=0\}$, $F_1(x,y)=(1,x a(x))$ and $F_2(x,y)=(0,ye^{\pphi(x,y)})$. Reducing $U$, if necessary, let $G_1(x,y)=(a_1(x,y),a_2(x,y)), G_2(x,y)=(b_1(x,y),b_2(x,y))$ be a frame 
for $\bD^2|_U$. Since $\bD^1(q)=\bD^2(q)$ at every point, we have $a_2(0,0)=b_2(0,0)=0$. Since $\bD^2(0,0)$ is one-dimensional, we may assume $a_1(0,0)\neq 0$. After possibly further reducing $U$, $\bD^2|_U=\mathrm{span}_{\con^\infty(M)}\{(1/a_1)G_1,G_2-(b_1/a_1)G_1\}$ and we  can assume 
$a_1(x,y)\equiv 1$ and  $b_1(x,y)\equiv 0$. Moreover, by $\bD^1(x,0)=\bD^2(x,0)$ we get  $a_2(x,0)=x a(x)$ and $b_2(x,0)\equiv 0$, whence 
 $a_2(x,y)=x a(x)+y \overline{a}_2(x,y)$ and $b_2(x,y)=y \overline{b}_2(x,y)$, with $\overline{a}_2,\overline{b}_2$ smooth functions.  Computing the Lie brackets  we get
$$
[G_1,G_2]|_{(x,0)}=(0,  x a  \overline{b}_2)|_{(x,0)},~~~[G_1,[G_1,G_2]]|_{(0,0)}=(0, a \overline{b}_2))|_{(0,0)}.
$$ 
Applying hypothesis \HH\ to ${\cal S}_2$ we have  $\overline{b}_2(x,0)\neq 0$ for all $x$ in a neighborhood of $0$.  Hence, up to reducing $U$,
\begin{eqnarray}
\bD^2|_U&=&\mathrm{span}_{\con^\infty(M)}\{G_1-(\overline{a}_2(x,y)/\overline{b}_2(x,y))G_2,(e^{\xi(x,y)}/\overline{b}_2(x,y))G_2\}
\nonumber\\
&=&\mathrm{span}_{\con^\infty(M)}\{F_1,F_2\}=\bD^1|_U.\nonumber
\end{eqnarray}
\hfill$\blacksquare$

{\bf Step 5.}
Thanks to Lemma \ref{coin} and  Proposition \ref{stessomodulo} we can assume $f_1=f_2=f$. In other words,  we reduce to the case ${\cal S}_1=(E,f,\langle\cdot,\cdot\rangle_1)$ and ${\cal S}_2=(E,f,\langle\cdot,\cdot\rangle_2)$. By compactness of $M$, 
there exists a constant $k\geq 1$ such that 
\begin{equation}\label{compat}
\frac{1}{k}\langle u,u\rangle_2\leq \langle u,u\rangle_1\leq k \langle u,u\rangle_2,\,\,\forall\,u\in E.
\end{equation}
For every $q\in M$ and $v\in\bD(q)$ let ${\bf G}^i_q(v)=\inf\{\langle u, u\rangle_i \mid u\in E_q,\f(u)=v\}$ (see section \ref{basdef}).
Clearly, 
\begin{equation}\label{confronto}
\frac{1}{k}{\bf G}^2_q(v)\leq {\bf G}^1_q(v)\leq k {\bf G}^2_q(v),\,\,\forall\,v\in f(E_q).
\end{equation}
By \r{compat}, 
admissible curves for ${\cal S}_1$ and ${\cal S}_2$ coincide. Moreover, given an admissible curve $\gamma:[0,T]\rightarrow M$, we can compare its length with respect to  ${\cal S}_1$ and ${\cal S}_2$ using \r{confronto}. Namely,
$$
\frac{1}{\sqrt{k}}\int_0^T\sqrt{{\bf G}^2_{\gamma(s)}(\dot\gamma(s))}ds\leq \int_0^T\sqrt{{\bf G}^1_{\gamma(s)}(\dot\gamma(s))}ds\leq 
\sqrt{k} \int_0^T\sqrt{{\bf G}^2_{\gamma(s)}(\dot \gamma(s))}ds.
$$
Since the Carnot-Caratheodory distance between two points is defined as the infimum of the lengths of the admissible curves joining them, 
we get 
$$
\frac{1}{\sqrt{k}}d_2(p,q)\leq d_1(p,q)\leq \sqrt{k} d_2(p,q),\,\,\forall\,p,q\in M.
$$
This is equivalent to say that the identity map is a Lipschitz equivalence between ${\cal S}_1$ and ${\cal S}_2$.
\hfill$\blacksquare$

\vspace{2cm}

\noindent{\bf Acnowledgements.} The authors are grateful to Andrei Agrachev for very helpful discussions.

\bibliographystyle{abbrv}
\bibliography{biblio_lipeq}

\end{document}